\documentclass[prl, sd, amsmath, amssymb, reprint]{revtex4-1}

\usepackage{amsmath}
\usepackage{amsfonts}
\usepackage{amssymb}
\usepackage{amsthm}
\usepackage{fancyhdr}
\usepackage{enumerate}
\usepackage[top=1in, bottom=1in, left=1in, right=1in]{geometry}
\usepackage{graphicx}
\usepackage[labelfont=bf]{caption}

\usepackage{caption}
\usepackage{subcaption}

\usepackage{verbatim}

\usepackage{wasysym}
\usepackage{hyperref}

\usepackage{wrapfig}

\begin{document}

\title{Volcano transition in a solvable model of oscillator glass}
\author{Bertrand Ottino-L\"{o}ffler and Steven H. Strogatz}
\affiliation{Center for Applied Mathematics, Cornell University, Ithaca, New York 14853}
\date{\today}

\begin{abstract}
In 1992 a puzzling transition was discovered in simulations of randomly coupled limit-cycle oscillators. This so-called volcano transition has resisted analysis ever since. It was originally conjectured to mark the emergence of an oscillator glass, but here we show it need not. We introduce and solve a simpler model with a qualitatively identical volcano transition and find, unexpectedly, that its supercritical state is not glassy. We discuss the implications for the original model and suggest experimental systems in which a volcano transition and oscillator glass may appear.
\end{abstract}

\maketitle

\section{Introduction}


Large systems of attractively coupled limit-cycle oscillators can show synchronization transitions analogous to ferromagnetic phase transitions~\cite{winfree1967biological, kuramoto1984chemical}. These transitions have been observed in chemical systems~\cite{kiss2002emerging} and are predicted for arrays of lasers~\cite{fabiny1993coherence, kozyreff2000laser, oliva2001laser}, biological oscillators~\cite{winfree1967biological}, Josephson junctions~\cite{wiesenfeld1996synchronization}, and optomechanical systems~\cite{heinrich2011optomechanical}. The analogy to ferromagnetism led Daido~\cite{daido1987population, daido1992quasientrainment} to conjecture that if the purely attractive  couplings were replaced by a frustrated mix of attractive and repulsive couplings, oscillator arrays could potentially behave like spin glasses~\cite{daido1987population, daido1992quasientrainment, fischer1993book, castellani2005spin}. So far, however, only a few counterparts of the phenomena observed in spin glasses have been seen in oscillator arrays~\cite{iatsenko2014glassy}. Finding, characterizing, and even defining a true ``oscillator glass'' remains controversial~\cite{daido1987population, daido1992quasientrainment, stiller1998dynamics, daido2000algebraic, stiller2000self, zanette2005synchronization, hong2011kuramoto, hong2011conformists, hong2012mean, iatsenko2014glassy}. 

The search for oscillator glass began with a natural model: a system of $N \gg 1$ phase oscillators with random couplings. The coupling strengths were chosen to be symmetric and Gaussian as in the Kirkpatrick-Sherrington spin-glass model~\cite{kirkpatrick1975solvable}. Simulations revealed that as the variance of the Gaussian couplings was increased, the model displayed a ``volcano'' transition~\cite{daido1992quasientrainment}. The name came from the shape of the model's two-dimensional, circularly symmetric distribution of complex local fields, which switched from being concave down at the origin to concave up, thus forming a volcano-like surface. Daido~\cite{daido1992quasientrainment} suggested this transition might signal the onset of an oscillator glass. Further evidence was provided by the numerical observation of slow (algebraic rather than exponential) relaxation from an initially synchronous state to an incoherent state. In an effort to explain these results analytically, later studies sought similar phenomena in more tractable models ~\cite{hong2011kuramoto, hong2011conformists, hong2012mean, iatsenko2014glassy, kloumann2014phase}, but so far the volcano transition and the glassy state have remained elusive. 

In this Letter we present a model with an exactly solvable volcano transition. It uses a coupling matrix whose rank is controlled by a parameter $K$. In the low-rank regime $2 \le K \ll  \log_2 N$ the model's dynamics, somewhat surprisingly, are non-glassy above threshold. Thus the volcano transition is not indicative of an oscillator glass; in the model studied here it merely signals a  synchronization transition in the presence of frustration. Unfortunately our analysis does not extend to the high-rank regime $K = O(N)$ of more direct relevance to Daido's results~\cite{daido1992quasientrainment}. For now that case remains out of reach. Whether a true oscillator glass exists in this or some other regime thus remains an open theoretical question.

Following Daido~\cite{daido1992quasientrainment}, our model consists of coupled phase oscillators. Oscillator $j$ couples with strength $J_{jk}$ to oscillator $k$ via the sine of their phase difference. The governing equations are
\begin{equation}\label{finiteN}
\dot\theta_j = \omega_j + \sum_{k = 1}^N J_{jk}\sin\left(\theta_k - \theta_j\right)
\end{equation}
for  $j = 1, \ldots, N$. Here $\theta_j$ denotes the phase of oscillator $j$ and $\omega_j$ is its natural frequency, selected at random from a given probability distribution. Instead of the Gaussian frequencies and couplings studied in Ref.~\cite{daido1992quasientrainment}, for the sake of solvability we consider Lorentzian-distributed frequencies with density 
\begin{equation*}\label{cauchy}
g(\omega) = \frac{1}{\pi(1+\omega^2)}
\end{equation*}
and define the couplings as follows. Given an even integer $K > 0$ and a coupling scale factor $J \ge 0$, let 
\begin{equation}\label{finiteJ}
J_{jk} = \frac{J}{N} \sum_{m=1}^K (-1)^m u^{(j)}_m u^{(k)}_m.
\end{equation}
Here, for each oscillator $j$ the interaction vector $(u_1^{(j)}, \ldots, u_K^{(j)})$ is a random binary vector of length $K$ with each entry independently being $\pm1$ with equal probability. Notice that the diagonal of $J_{jk}$ will always be 0 (since it is an alternating sum of 1's), and $J_{jk} = J_{kj}$. In the limit of large $N$ the parameter $K$ equals the rank of the coupling matrix $J_{jk}$. Furthermore, if we fix $K \equiv N$ and let $N$ get large, the off-diagonal entries converge to normal random variables with a standard deviation of $J/\sqrt{N}$. So when $K = N \gg 1$, our  construction approximates Daido's original Gaussian couplings.

\begin{figure}
\centering
\includegraphics[width = 0.5\textwidth]{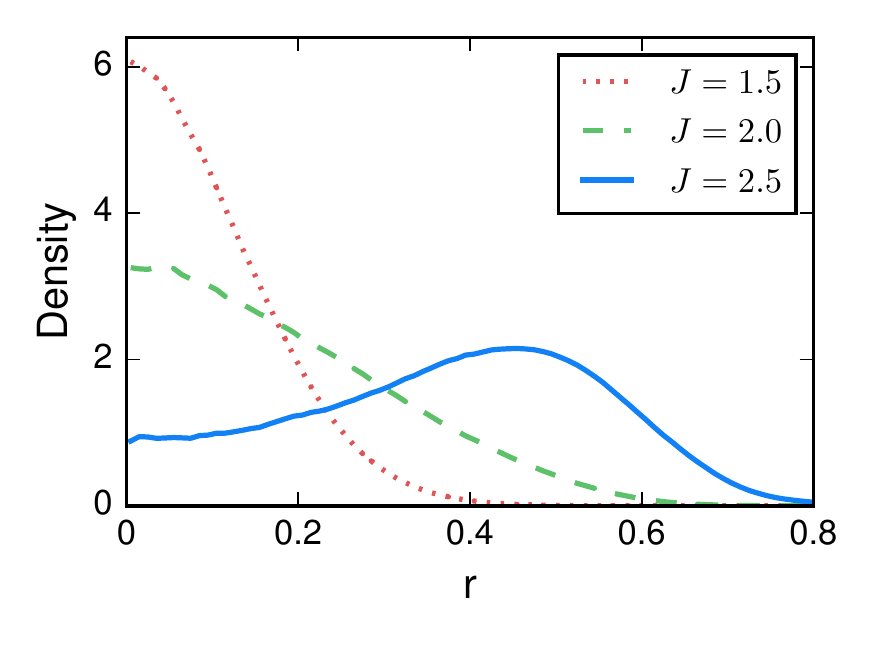}
\caption{Radial distribution of local fields. Each curve represents the averaged density over 500 simulations of Eq.~\eqref{finiteN}, using $N = 250$, $K = 4$, fourth-order Runge-Kutta integration with a step size of 0.01, 1000 transient steps, 2000 recorded steps, and uniformly random initial phases.}
\label{Plot_PDists}
\end{figure}

\begin{figure}
\centering
\includegraphics[width = 0.5\textwidth]{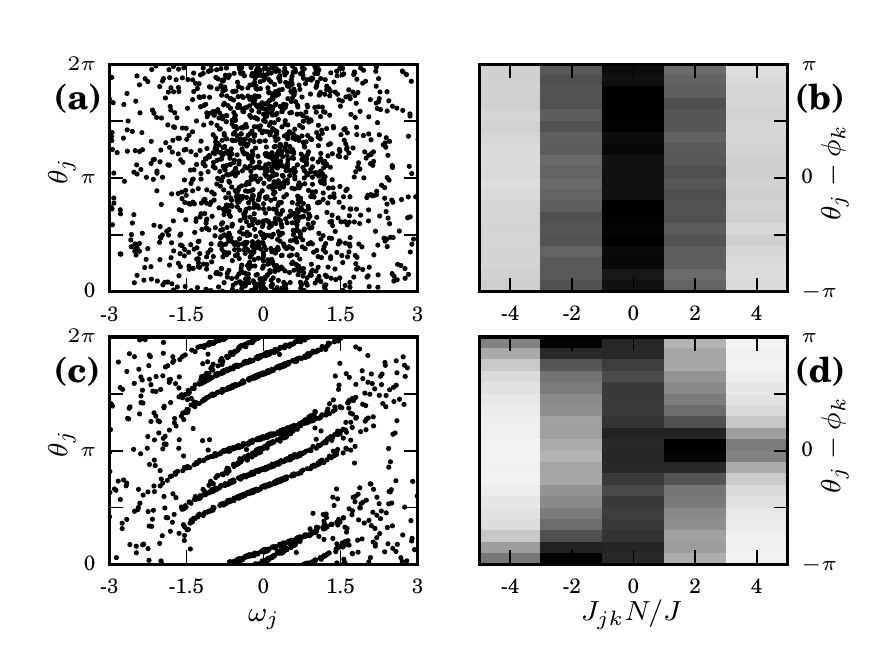}
\caption{Oscillator phase distributions below and above the volcano transition. In (a) and (b), $J=1$; in (c) and (d), $J=3$. Each panel shows results for simulations of $N = 2000$ and $K = 6$; other parameters as in Fig.~\ref{Plot_PDists}.  (a) Below the volcano transition, the system is incoherent. (b) Density of $\theta_j - \phi_k$, indicating an oscillator's phase relative to that of the local field angle, plotted against the associated coupling strength $J_{jk}$, normalized and averaged across all $k$. Darker shades represent higher density. The uniform vertical stripes show that when $J=1$ the local field has negligible influence on oscillator phases. (c)  Above the volcano transition, phase-locked clusters appear. (d) Dark horizontal bands at $\theta_j - \phi_k = 0$ and $\pm \pi$ indicate tendency of oscillators to align or anti-align to local field phases, depending on the sign of $J_{jk}$. }
\label{Plot_Thetas}
\end{figure}

To show numerically that our model has a volcano transition, we compute its complex local fields~\cite{daido1992quasientrainment} 
\begin{equation*}\label{localField}
P_j = r_j e^{i\phi_j} := \sum_{k = 1}^N J_{jk} e^{i\theta_k},
\end{equation*}
for $j = 1, \ldots, N$. Equation~\eqref{finiteN} then becomes
\begin{equation*}
\dot\theta_j = \omega_j + r_j \sin\left(\phi_j - \theta_j\right).
\end{equation*}
By keeping track of the $P_j$ over time, we obtain a distribution of their magnitudes $r_j$ for each realization of $\omega$ and $J_{jk}$. Figure~\ref{Plot_PDists} averages these distributions over many realizations. As $J$ increases from 1.5 to 2.5, the distribution changes from concave down at the origin to concave up and volcano-like. At a critical $J_c$, the origin no longer attracts the maximum density. This $J_c$ defines the volcano transition.

Figure~\ref{Plot_Thetas} illustrates how the individual oscillator phases $\theta_j$ behave on either side of the transition. For $J < J_c$ the system is incoherent [Fig.~\ref{Plot_Thetas}(a)]. The phases of the oscillators are uniformly distributed and bear no relation to the coupling strength $J_{jk}$ or the phase $\phi_k$ of the complex local field [Fig.~\ref{Plot_Thetas}(b)]. In contrast, for $J > J_c$ the oscillators with small $|\omega_j|$ form phase-coherent clusters [Fig.~\ref{Plot_Thetas}(c)]. Figure~\ref{Plot_Thetas}(d) suggests that this partial synchronization is induced by the local fields: if oscillator $j$ couples positively (attractively) to oscillator $k$, then oscillator $j$ tends to align with the $k$th local field, whereas if they are negatively (repulsively) coupled, then oscillator $j$ tends to anti-align with the local field. In some realizations we have also observed clustering at phase differences other than 0 and $\pi$, for moderate values of $J_{jk}$.

Turning now to the analytical results, we examine Eqs.~\eqref{finiteN} and \eqref{finiteJ} in the continuum limit $N \rightarrow \infty$ with $K$ held fixed. Using an Eulerian description, we replace our discrete system of oscillators with a continuous fluid moving around the unit circle. Its state is described by a density $f(\theta, \omega, u, t)$ of oscillators with phase $\theta$, natural frequency $\omega$, and interaction vector $u$. In this framework the dynamics are given by a continuity equation $f_t + {(f \nu)}_{\theta} = 0$, where the subscripts denote partial differentiation, $\nu$ represents the velocity field on the circle given by the continuum limit of Eq.~\eqref{finiteN},
\begin{equation}\label{infiniteN} 
\nu(\theta, \omega, u, t) = \omega + \langle J(u, u') \sin(\theta' - \theta) \rangle,
\end{equation}
and $\langle\cdot\rangle$ denotes integration using the time-dependent measure $ f(\theta', \omega', u', t) d\theta' g(\omega')d\omega' \rho(u')du'$. The coupling term $J(u,u')$  in Eq.~\eqref{infiniteN} plays the role of $J_{jk}$ in Eq.~\eqref{finiteJ}. It is given by 
\begin{equation*}\label{infiniteJ}
J(u,u') := J \sum_{m=1}^K (-1)^m u_{m} u_{m}'.
\end{equation*}
As before, $u$ and $u'$ are random interaction vectors of length $K$, all of whose entries are $\pm 1$ with probability 1/2 each. Thus the probability of any particular vector is $2^{-K}$. The associated term in the measure is $\rho(u') = 2^{-K} \sum_{v} \delta(u'-v)$, where the sum runs over all the equally likely $v \in \{\pm 1\}^K$. Similarly, the continuum limit of the local field is  
\begin{equation*}\label{orderParam}
P(u,t) = \langle J(u,u')e^{i\theta'} \rangle.
\end{equation*}
Inserting $P$ in Eq.~\eqref{infiniteN} gives 
\begin{equation*}\label{infiniteN2}
\nu(\theta, \omega, u, t) = \omega + \frac{1}{2i}\left[ e^{-i\theta} P(u, t) - \mbox{c.c.} \right].
\end{equation*}

Having derived the continuum model, we reduce it with the Ott-Antonsen ansatz~\cite{ott2008low, ott2009long, pikovsky2008partially, pikovsky2015review}, a technique that yields  the exact long-term dynamics of Kuramoto oscillator models with sinusoidal coupling and Lorentzian frequencies. Following the standard procedure we seek solutions of the form 
\begin{equation*}\label{ottAnsatz}
f(\theta, \omega, u, t) = \frac{1}{2\pi} \left[1 + \sum_{n=1}^\infty \alpha(\omega, u, t)^n e^{int} + \mbox{c.c.} \right]
\end{equation*}
and define $a(u,t) := \alpha(-i, u, t)$. Then we find
\begin{equation}\label{ottODE}
\dot a(u,t) = -a(u,t) + \frac{ P^*(u,t) - a(u,t)^2P(u,t)}{2}
\end{equation}
where
\begin{alignat*}{1}\label{PIntegral}
P(u,t) &= \int J(u,u') a^*(u',t)\rho(u')du' \\
& = \frac{J}{2^K} \sum_{u'} \sum_{m=1}^K (-1)^m u_m u'_m a^*(u', t).
\end{alignat*}
Finally, by replacing $P$ in Eq.~\eqref{ottODE} with this sum, we get a closed set of $2^K$ ordinary differential equations for the $a(u,t)$, one for each possible choice of $u$. 

Equation~\eqref{ottODE} has rich dynamics, but for our purposes it suffices to analyze the stability of its trivial fixed point, $a(u,t)=0$ for all $u$ and $t$, because this state corresponds to the incoherent state of Eq.~\eqref{finiteN}. The volcano transition occurs precisely when this state goes unstable. Thus, to calculate $J_c$ we linearize Eq.~\eqref{ottODE} about $a \equiv 0$ and determine when one of its eigenvalues is 0. 
The Jacobian is 
\begin{equation}\label{Jacobian}
-I + \frac{J}{2^{K+1}} A. 
\end{equation}
Here $I$ is the $2^K \times 2^K$ identity matrix and 
\begin{equation*}\label{JacMat}
A_{u,v} = \sum_{m=1}^K (-1)^m u_{m} v_{m} 
\end{equation*} 
where the entries of $A$ have been conveniently indexed by binary strings $u, v \in \{\pm 1\}^K$. 
The eigenvalues of $A$ can be found explicitly. To do so we write down the eigenvectors (which we guessed by generalizing from small examples) and then read off the eigenvalues. For each integer $1 \le n \le  K$ and each binary string $v \in \{\pm 1\}^K$, define a vector $\zeta^{(n)} \in \mathbb{R}^{2^K}$ whose $v$th entry is $\zeta^{(n)}_v = v_{n}$.
One can check that the set of all $K$ such vectors is orthogonal and, by using the evenness of $K$, that $A \zeta^{(n)} = (-1)^n 2^K \zeta^{(n)}$.
Moreover, given any $\eta$ perpendicular to all the $\zeta^{(n)}$, one finds $A \eta = 0$.
Therefore, $A$ has exactly three distinct eigenvalues: $+2^K$ with multiplicity $K/2$, $-2^K$ with multiplicity $K/2$, and $0$ with multiplicity $2^K - K$. 
Consequently the Jacobian~\eqref{Jacobian} has three distinct eigenvalues, with the largest always being $-1+J/2$. The conclusion is that the incoherent state for the continuum model loses stability at 
\begin{equation}\label{Critical}
J_c = 2.
\end{equation}
This result holds for \emph{any} even value of $K$. 

The next question is whether $J_c=2$ gives a good approximation to $J_c$ when $N$ is finite. To anticipate the answer, recall that the continuum model reduces to the $2^K$-dimensional system~\eqref{ottODE}. For the finite-$N$ system \eqref{finiteN} to have any chance of behaving like a continuous fluid of oscillators, we need it to have many oscillators per $u$, and hence to have $N \gg 2^K$. 

To test these ideas we simulate the finite-$N$ system and estimate $J_c$ carefully. To pinpoint the volcano transition we first compute the one-dimensional (1D) distribution of local field magnitudes $r_j \ge 0$ and fit it to the sum of two normal distributions, with one centered at $\mu$ and the other at $-\mu$, and both with variance $\sigma^2$. In other words, we approximate the 1D density of local field magnitudes by
\begin{equation*}\label{distAnsatz2}
h(r) = \frac{2}{\sqrt{2\pi\sigma^2}}\exp\left(\frac{-\mu^2 - r^2}{2\sigma^2}\right) \cosh\left(\frac{\mu r}{\sigma^2}\right),
\end{equation*}
for $r \ge 0$. To obtain the full 2D distribution of the $P_j$'s, we impose azimuthal symmetry by rotating and rescaling the 1D density above. 

\begin{figure}
\centering
\includegraphics[width = 0.5\textwidth]{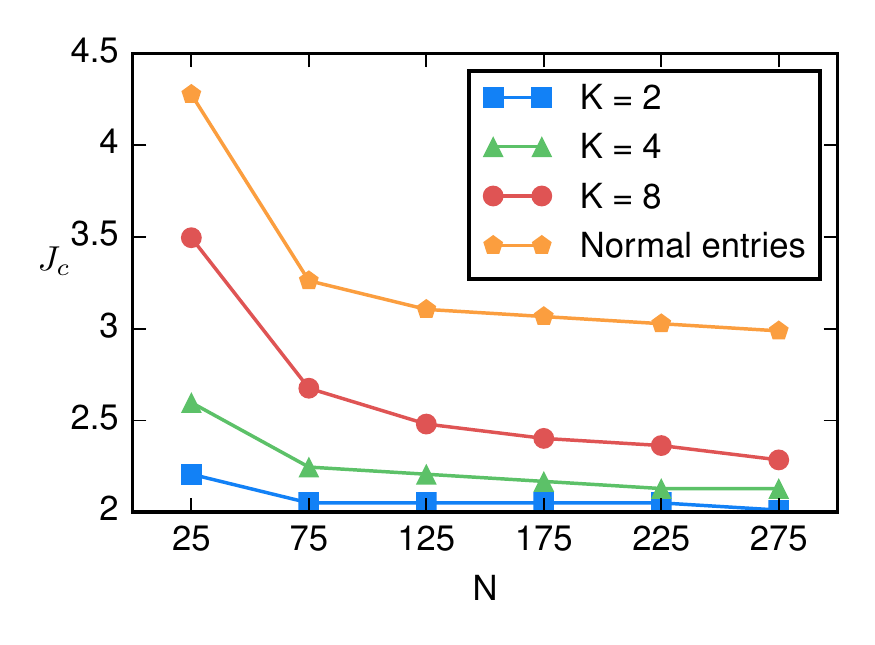}
\caption{Critical value $J_c$ versus $N$ and $K$. Each value of $J_c$ was estimated by using a bisection method on the value of $M_{+1}M_{-1}$ to achieve an accuracy of $\lesssim 0.02$. For each $J$ we sample $J_{jk}$ at least 100 times, simulate Eq.~\eqref{finiteN}, evaluate $M_{+1}M_{-1}$, and keep track of the running standard deviation of these products. If the current value of $M_{+1}M_{-1}$ is more than 1.5 standard deviations from 1.4694, the bisection continues; otherwise further simulations are run, up to a maximum of $10^5$ simulations. Each simulation consists of 1000 transient steps followed by 2000 recorded steps of a fourth-order Runge Kutta integration with a step size of 0.01, with initial phases all set to 0.}
\label{Plot_CritVrsN}
\end{figure}

\begin{figure}
\centering
\includegraphics[width = 0.5\textwidth]{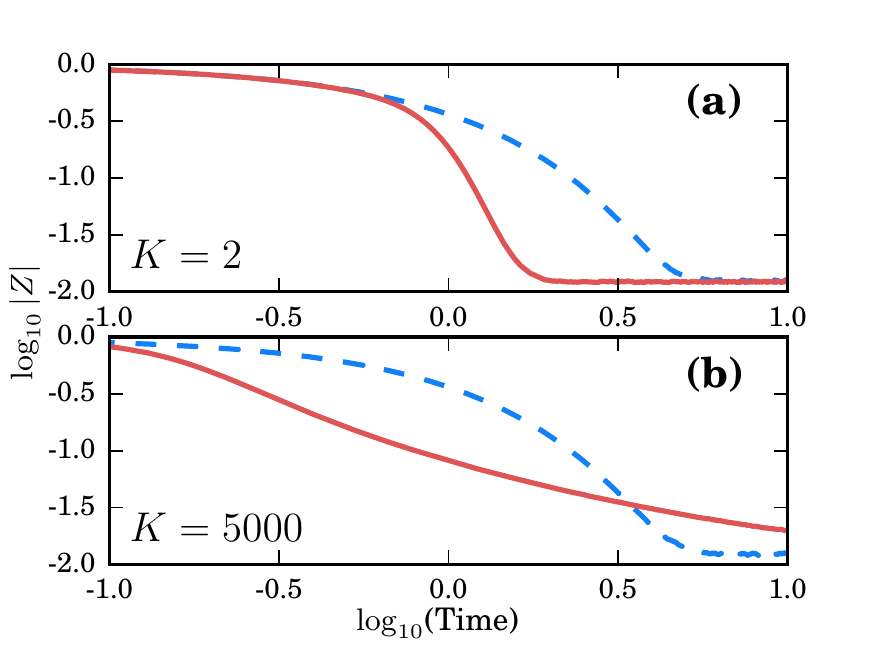}
\caption{Log-log plot for the decay of the order parameter $Z(t)$. Each curve is the average of 750 numerical integrations of Eq.~\eqref{finiteN} for  $N = 5000$ oscillators starting from the in-phase state ($\theta_j=0$ for all $j$) and run for 1000 steps with a step size of 0.01. Solid curves show coupled systems with $J = 10$; dashed curves show uncoupled systems with $J= 0$ for which the order parameter decays exponentially: $Z(t) = e^{-t}$. (a) Low-rank regime: $K \ll \log_2(N)$. For 
$K = 2$,  $Z(t)$ decays exponentially down to the noise floor. Exponential decay is expected in this regime because the dynamics of Eq.~\eqref{finiteN} are well approximated by the low-dimensional  system \eqref{ottODE}. (b) High-rank regime: $K=N=5000$. When $K = O(N)$ and $J > J_c$, the relaxation of $Z$ slows markedly, resembling the algebraic decay in glass.
}
\label{Plot_Decay}
\end{figure}


The functional form of $h(r)$ allows us to identify its convexity at the origin easily. It is concave down when $\gamma:= \mu^2/\sigma^2 < 1$ and concave up when $\gamma > 1$. 
To measure $\gamma$ numerically, we use the method of moments on the 2D distribution and find that the product of the first and negative first moments is 
\begin{equation*}\label{momentFit}
M_{+1}M_{-1} = \frac{\pi}{2} \frac{1+\gamma}{[e^{-\gamma/2} + \sqrt{\pi\gamma/2} \mbox{Erf}(\sqrt{\gamma/2}) ]^2}.
\end{equation*}
The left hand side can be numerically estimated by aggregating moments from multiple simulations, along with an appropriate estimate of an error on its total. The right hand side can be proven to imply that $h(r)$ is concave down at the origin (and therefore $J < J_c$) if and only if $M_{+1}M_{-1} \gtrsim 1.4694$. Thus by measuring these two moments we can use a bisection algorithm to zero in on $J_c$. 

Figure~\ref{Plot_CritVrsN} shows that when $K$ is small, $J_c = 2$ becomes an increasingly good estimate as $N$ gets large. For comparison we also computed $J_c$ for a Gaussian coupling model in which $J_{jk}$ is a random symmetric matrix with normally distributed entries having mean zero and variance $J^2/N$. As noted earlier, our coupling matrix \eqref{finiteJ} converges to such a Gaussian matrix when $K=N \gg 1$, but our analytical approach does not extend to this large-$K$ regime. So although the value of $J_c$ for Gaussian coupling decreases as $N$ gets large, we cannot predict whether $J_c$ asymptotically approaches 2 or not. 

A proposed signature feature of oscillator glasses~\cite{daido1992quasientrainment, iatsenko2014glassy, daido2000algebraic, stiller2000self, stiller1998dynamics, iatsenko2014glassy} is nonexponential relaxation of the order parameter 
\begin{equation*}
Z(t) := \sum_{k = 1}^N e^{i\theta_k(t)} .
\end{equation*}
Figure~\ref{Plot_Decay} plots the decay of the order parameter for our model. In the low-rank regime $K \ll \log_2 N$ to which our continuum theory applies, Fig.~\ref{Plot_Decay}(a) shows that $Z$ decays exponentially fast. This is to be expected, given that the dynamics reduce to a low-dimensional set of ordinary differential equations~\eqref{ottODE} in this regime. So the dynamics are not glassy here, even  above the volcano transition. 
However, outside the low-rank regime there is some indication that the model can exhibit a glassy state. Figure~\ref{Plot_Decay}(b) shows that when $K=N \gg 1$, the order parameter $Z(t)$ decays roughly algebraically for sufficiently large $J$. This finding is consistent with results from the closely related Gaussian coupling model, which has been claimed (controversially) to have algebraic decay~\cite{daido1992quasientrainment, daido2000algebraic, stiller2000self, stiller1998dynamics}. Analytically understanding the nature of this decay, in both our model and Daido's, remains an open problem. Another important future direction is the experimental investigation of oscillator glasses. The most promising experimental setup in which to search for them may be a large array of photosensitive chemical oscillators coupled through a programmable spatial light modulator, as recently used \cite{totz2017spiral} to demonstrate the  existence of spiral wave chimeras. 

We thank Hiroaki Daido for helpful interactions. This research was supported by a Sloan Fellowship to Bertrand Ottino-L\"{o}ffler in the Center for Applied Mathematics at Cornell, as well as by NSF grant DMS-1513179 to Steven Strogatz.

\bibliography{OscillatorGlass_Bib_3}{}

\end{document}